\documentclass[]{article}
\usepackage{latexsym}
\usepackage{amssymb}
\parskip=\medskipamount
\newtheorem{definition}{Definition}[section]
\newtheorem{lemma}[definition]{Lemma}
\newtheorem{theorem}[definition]{Theorem}
\newtheorem{proposition}[definition]{Proposition}
\newtheorem{corollary}[definition]{Corollary}
\newtheorem{remark}[definition]{Remark}
\newtheorem{example}[definition]{Example}

 1   
 1  

\def\QED{\hskip0.1em\hfill\null\ \null\nobreak\hfill
\kern3pt\lower1.8pt\vbox{\hrule\hbox
{\vrule\kern1pt\vbox{\kern1.7pt \hbox{$\scriptstyle
QED$}\kern0.2pt}\kern1pt\vrule}\hrule}}

\def\lcf{\lbrack\! \lbrack}
\def\rcf{\rbrack\! \rbrack}

\mathchardef\za="710B  
\mathchardef\zb="710C  
\mathchardef\zg="710D  
\mathchardef\zd="710E  
\mathchardef\zve="710F 
\mathchardef\zz="7110  
\mathchardef\zh="7111  
\mathchardef\zvy="7112 
\mathchardef\zi="7113  
\mathchardef\zk="7114  
\mathchardef\zl="7115  
\mathchardef\zm="7116  
\mathchardef\zn="7117  
\mathchardef\zx="7118  
\mathchardef\zp="7119  
\mathchardef\zr="711A  
\mathchardef\zs="711B  
\mathchardef\zt="711C  
\mathchardef\zu="711D  
\mathchardef\zvf="711E 
\mathchardef\zq="711F  
\mathchardef\zc="7120  
\mathchardef\zw="7121  
\mathchardef\ze="7122  
\mathchardef\zy="7123  
\mathchardef\zf="7124  
\mathchardef\zvr="7125 
\mathchardef\zvs="7126 
\mathchardef\zf="7127  
\mathchardef\zG="7000  
\mathchardef\zD="7001  
\mathchardef\zY="7002  
\mathchardef\zL="7003  
\mathchardef\zX="7004  
\mathchardef\zP="7005  
\mathchardef\zS="7006  
\mathchardef\zU="7007  
\mathchardef\zF="7008  
\mathchardef\zW="700A  

\newcommand{\be}{\begin{equation}}
\newcommand{\ee}{\end{equation}}
\newcommand{\ra}{\rightarrow}

\newcommand{\bea}{\begin{eqnarray}}
\newcommand{\eea}{\end{eqnarray}}
\newcommand{\beas}{\begin{eqnarray*}}
\newcommand{\eeas}{\end{eqnarray*}}
\newcommand{\Z}{{\mathbb Z}}
\newcommand{\R}{{\mathbb R}}

\newcommand{\we}{\wedge}
\newcommand{\nn}{\nonumber}
\newcommand{\ot}{\otimes}

\newcommand{\pa}{\partial}
\newcommand{\ti}{\times}

\newcommand{\Li}{{\cal L}}

\def\la{\langle}
\def\ran{\rangle}
\def\wt{\widetilde}

\begin{document}

\title{Poisson-Jacobi reduction of homogeneous tensors
\thanks{Research supported by the Polish Ministry of Scientific Research
and Information Technology under the grant No. 2 P03A 020 24 and MCYT
grants BFM2000-0808 and BFM2003-01319. D.~Iglesias wishes to thank the
Spanish Ministry of Education and Culture and Fulbright program for a FPU
grant and for a MECD/Fulbright postdoctoral grant. }}

\author{J. Grabowski$^{1}$, D. Iglesias$^{2}$,  J.C.
Marrero$^{3},$ E. Padr\'on$^{3}, $ P. Urbanski$^{4}$
\\[10pt]
 {\small\it$^1$ Mathematical Institute, Polish Academy of
Sciences}\\[-6pt] {\small\it \'Sniadeckich 8, P.O.Box 21, 00-956
Warsaw, Poland }\\[-6pt]{\small\it E-mail:
jagrab@impan.gov.pl}\\[-4pt]{\small\it $^2$Department of
Mathematics,  The Pennsylvania State University} \\[-6pt]
{\small\it University Park, PA 16802, USA}\\[-6pt] {\small\it
E-mail: iglesias@math.psu.edu}\\[-4pt] {\small\it $^3$Departamento de
Matem\'atica Fundamental, Facultad de Matem\'aticas}\\[-6pt]
{\small\it Universidad de la Laguna, La Laguna} \\[-6pt]
{\small\it Tenerife, Canary Islands, SPAIN}\\[-6pt] {\small\it
E-mail: jcmarrer@ull.es, mepadron@ull.es
}\\[-4pt] {\small\it$^4$Division of Mathematical Methods in
Physics, University of  Warsaw}\\[-6pt] {\small\it Ho\.za 74, 00-682
Warsaw, Poland }\\[-6pt] {\small\it E-mail: urbanski@fuw.edu.pl} }
\date{\empty}

\maketitle
\begin{abstract}
\noindent The notion of homogeneous tensors is discussed. We show that
there is a one-to-one correspondence between multivector fields on a
manifold $M$, homogeneous with respect to a vector field $\Delta$ on $M$,
and first-order polydifferential operators on a closed submanifold $N$ of
codimension $1$ such that $\Delta$ is transversal to $N$. This
correspondence relates the Schouten-Nijenhuis bracket of multivector
fields on $M$ to the Schouten-Jacobi bracket of first-order
polydifferential operators on $N$ and generalizes the Poissonization of
Jacobi manifolds. Actually, it can be viewed as a super-Poissonization.
This procedure of passing from a homogeneous multivector field to a
first-order polydifferential operator can be also understood as a sort of
reduction; in the standard case -- a half of a Poisson reduction. A dual
version of the above correspondence yields in particular the
correspondence between $\Delta$-homogeneous symplectic structures on $M$
and contact structures on $N$.
\end{abstract}

\begin{quote}

\vspace{-5pt}

\end{quote}

\setcounter{section}{0}

\section{Introduction}

As it has been observed in \cite{KoS}, a Lie algebroid structure
on a vector bundle $E$ can be identified with a Gerstenhaber
algebra structure on the exterior algebra of multisections of $E$,
$Sec(\wedge E),$ which is just a graded Poisson bracket (Schouten
bracket) on $Sec(\wedge E)$ of degree $-1$, that is, the Schouten
bracket is graded commutative, satisfies the graded Jacobi
identity and the graded Leibniz rule.

In the particular case of the Lie algebroid structure on the
tangent vector bundle of an arbitrary manifold $M$  one obtains
the Schouten-Nijenhuis bracket $\lcf\cdot,\cdot\rcf_M$ on the
space of multivectors on $M$.

For a graded commutative algebra with $1$, a natural
generalization of a graded Poisson bracket is a graded Jacobi
bracket: we replace the graded Leibniz rule by that $\{a,\cdot\}$
is a first-order differential operator on ${\cal A}$, for every
$a\in {\cal A}$ (cf. \cite{GM2}).

Graded Jacobi brackets on $Sec(\wedge E)$ of degree $-1$ are
called Schouten-Jacobi brackets. These brackets are in one-to-one
correspondence with pairs $(E,\phi_0)$, where $\phi_0\in Sec(E^*)$
is a $1$-cocycle in the Lie algebroid cohomology of $E$. In this
case, we said that $(E,\phi_0)$ is a generalized Lie algebroid
(Jacobi algebroid) (see \cite{GM,IM}).

A canonical example of a Jacobi algebroid is $(T^1M,(0,1))$ where
$T^1M=TM\oplus \R$ is the Lie algebroid of first-order
differential operators on the space of smooth functions on $M,$
$C^\infty(M),$ with the bracket $$\lcf X\oplus f,Y\oplus g\rcf_M^1
=[X,Y]\oplus (X(g)-Y(f)),$$ for $X\oplus f,Y\oplus g\in Sec(T^1M)$
(see \cite{M,N}) and the $1$-cocycle $\phi_0=(0,1)\in
\Omega^1(M)\oplus C^\infty(M)$.

It is well-known that a Poisson structure on a manifold $M$ can be
interpreted as a canonical structure for the Schouten-Nijenhuis
bracket $\lcf\cdot,\cdot\rcf_M$ of multivector fields on $M$,
i.e., as an element $\Lambda\in Sec(\wedge^2TM)$ satisfying the
equation $\lcf\Lambda,\Lambda\rcf_{M}=0.$ In similar way, a Jacobi
structure is a canonical structure for the Jacobi bracket
$\lcf\cdot,\cdot\rcf^1_M$.

On the other hand, it is proved in \cite{DLM} that if $\Lambda$ is a
homogeneous Poisson tensor with respect to a vector field $\Delta$ on the
manifold $M$ and $N$ is a $1$-codimensional closed submanifold of $M$ such
that $\Delta$ is transversal to $N$ then $\Lambda$ can be reduced to a
Jacobi structure on $N$.

The main purpose of this paper is to give an explicit (local)
correspondence between $\Delta$-homogeneous multivector fields on $M$ and
first-order polydifferential (i.e. skew-symmetric multidifferential)
operators on $N.$ This correspondence relates the Schouten-Nijenhuis
bracket of multivector fields on $M$ to the Schouten-Jacobi bracket of
first-order polydifferential operators on $N$. This is of course a
generalization of \cite{DLM} formulated in a structural way. It explains
the role of homogeneity for certain reduction procedures, e.g. in passing
from Poisson to Jacobi brackets (in mechanics: from symplectic form to a
contact form). But our result can be applied in Nambu-Poisson geometry
(cf. Corollary 3.13) or multisymplectic geometry and classical field
theories as well.

The paper is organized as follows. In Section 2 we recall the
notions of Schouten-Nijenhuis and Schouten-Jacobi brackets
associated with any smooth manifold. In Section 3.1 we introduce
the notion  of $\Delta$-homogeneous tensors on a homogeneous
structure $(M,\Delta)$ (a pair where $M$ is a manifold and
$\Delta$ is a vector field on $M$).

Moreover, for a particular class of homogeneous structures (strict
homogeneous structures), we will characterize  the
$\Delta$-homogeneous contravariant $k$-tensors in terms of their
corresponding $k$-ary brackets.

The main result of the paper is Theorem \ref{main} of Section 3.2, which
provides the one-to-one correspondence between homogeneous multivector
fields and polydifferential operators we have already mentioned. This
result is a generalization of the result of \cite{DLM} cited above and it
allows us also to relate homogeneous Nambu-Poisson tensors on $M$ to
Nambu-Jacobi tensors on $N$. These results are local. We obtain global
results in the particular case of the Liouville vector field
$\Delta=\Delta_E$ of a vector bundle $\tau:E\to M$. We called this
correspondence a Poisson-Jacobi reduction, since it can be understood as a
sort of reduction, a half of a Poisson reduction (cf. Remark 3.12, ii)).

Finally, we prove a dual version of Theorem 3.11. What we get is a
one-to-one correspondence between homogeneous differential forms on $M$
and elements of $Sec(\we(T^*N\oplus\R))$ represented by pairs
$(\alpha^0,\alpha^1)$, where $\alpha^0$ is a $k$-form on $N$ and
$\alpha^1$ is a $(k-1)$-form on $N$. This correspondence relates the de
Rham differential on $M$ with deformed Lie algebroid differential
associated with the Schouten-Jacobi bracket $\lcf\cdot,\cdot\rcf_{M}^1$
(see \cite{IM,GM}).

Note that the Grassmann algebra $Sec(\we TM)$ can be viewed as the
algebra of functions on the supermanifold $\zP T^*M$ (the space of
the cotangent bundle to $M$ with reversed parity of fibers, cf.
\cite{AKSZ}) the Schouten-Nijenhuis bracket on $Sec(\we TM)$
represents the canonical (super) Poisson bracket on $\zP T^*M$. In
this picture, the equation $\lcf\Lambda,\Lambda\rcf_{M}=0$ for a
Poisson tensor $\Lambda$ is just a particular case of the Master
Equation in Batalin-Vilkovisky formalism. The algebraic structure
of Batalin-Vilkovisky formalism in field theories (see \cite{Ge})
have been recognized as a homologic vector field generating a
Schouten-Nijenhuis-type bracket on the corresponding graded
commutative algebra like the Schouten-Nijenhuis bracket
(Gerstenhaber algebra) of a Lie algebroid \cite{KoS,KS2}. The
Schouten-Jacobi bracket can be regarded as a super-Jacobi bracket,
so Theorem 3.11 can be understood as a super or fermionic version
of the original result \cite{DLM}. Note also that higher-order
tensors represent higher-order operations on the ring of
functions. Together with the Schouten-Nijenhuis or Schouten-Jacobi
bracket, possibly for higher gradings, this can be a starting
point for certain strongly homotopy algebras (cf. the paper
\cite{St} by J.~Stasheff who realized that homotopy algebras
appear in string field theory). A relation of some strongly
homotopy algebras with Batalin-Vilkovisky formalism was discovered
by B.~Zwiebach and applied to string field theory \cite{Zw}.
Theorem 3.11 means that in homogeneous cases we can reduce the
structure to the same super Lie bracket on a smaller manifold. The
difference is that we deal not with derivations but with
first-order differential operators. The structure of the
associative product is deformed by this bracket isomorphism, so we
get not a super Poisson but a super Jacobi bracket. On the level
of differential forms this corresponds to a deformation of the de
Rham differential of the type $d^1\zm=d\zm+\phi\we\zm$, where
$\phi$ is a closed 1-form. This is exactly what was already
considered by E.~Witten \cite{Wi} and used in studying of spectra
of Laplace operators.

\section{Graded Lie brackets}\label{graded}

In this section we will recall several natural graded Lie brackets
of tensor fields associated with any smooth manifold $M$. First of
all, on the tangent bundle $TM$, we have a Lie algebroid bracket
$[\cdot ,\cdot ]$ defined on the space $\mathfrak X (M)$ of vector
fields -- derivations of the algebra $C^\infty(M)$ of smooth
functions on $M$.

If $A(M)=\oplus _{k\in \Z} A^k(M)$ is the space of multivector
fields (i.e., $A^k(M)=Sec (\wedge ^k TM)$) then we can define the
Schouten-Nijenhuis bracket (see \cite{Sc,Ni}) $\lcf \cdot ,\cdot
\rcf _M:A^p(M)\times A^q(M)\to A^{p+q-1}(M)$ as the unique graded
extension to $A(M)$ of the bracket $[\cdot ,\cdot ]$ of vector
fields, such that:

{\it i)} $\lcf X,f \rcf _M =X(f)$, for $X\in \mathfrak X (M)$ and
$f\in C^\infty (M)$;

{\it ii)} $\lcf P,Q\rcf _M =-(-1)^{(p-1)(q-1)}\lcf Q,P\rcf _M$, for $P\in
A^{p}(M)$, $Q\in A^{q}(M)$.

{\it iii)} $\lcf P,Q\wedge R\rcf _M =\lcf P,Q\rcf _M\wedge
R+(-1)^{(p-1)q}Q\wedge\lcf P,R\rcf _M$, for $P\in A^{p}(M)$, $Q\in
A^{q}(M)$ and $R\in A^{\ast }(M)$;

{\it iv)} $(-1)^{(p-1)(r-1)}\lcf\lcf P,Q\rcf _M,R\rcf
_M+(-1)^{(p-1)(q-1)}\lcf\lcf Q, R \rcf _M, P \rcf _M+(-1)^{(q-1)(r-1)}$
$\lcf \lcf R,P\rcf _M,$ $Q\rcf _M=0$, for $P\in A^{p}(M)$, $Q\in A^{q}(M)$
and $R\in A^{r}(M)$.

On the other hand, if $\Omega(M)=\oplus_{k\in \Z}\Omega^k(M)$ is
the space of differential forms (that is,
$\Omega^k(M)=Sec(\wedge^k(T^*M)))$ we can consider the usual
differential $d_M:\Omega^p(M)\to \Omega^{p+1}(M)$ as the map
characterized by the following properties:
\begin{enumerate}
\item $d_M$ is a $\R$-linear map.
\item $d_M(f)$ is the usual differential of $f$, for $f\in
C^\infty(M)$.
\item $d_M(\alpha\wedge \beta)=d_M\alpha\wedge \beta +
(-1)^p\alpha\wedge d_M\beta$, for $\alpha\in \Omega^p(M)$ and
$\beta\in \Omega^q(M)$.
\item $d_M^2=0,$ that is, $d_M$ is a cohomology operator.
\end{enumerate}

In a similar way, on the bundle of first-order differential
operators on $C^\infty (M)$, $T^1M=TM\oplus \R$, there exists a
Lie algebroid bracket given by

\begin{equation}\label{firstorder}
\lcf X\oplus f,Y\oplus g\rcf_M^1 =[X,Y]\oplus (X(g)-Y(f)),
\end{equation}
for $X\oplus f,Y\oplus g\in Sec(T^1M)$ (see \cite{M,N}).

The space $D^k(M)=Sec (\wedge ^k(T^1M))$ of sections of the vector bundle
$\wedge ^k (T^1M) \to M$ can be identified with $A^k(M)\oplus A^{k-1}(M)$
in the following way. If $I_M=0\oplus 1_M\in Sec (T^1M)$ and $\phi _M \in
Sec ((T^1M)^\ast )$ is the ``canonical closed 1-form'' defined by
$\phi_M(X\oplus f)=f$, then there exists an isomorphism between $D^k(M)$
and $A^k(M)\oplus A^{k-1}(M)$ given by the formula:

\[
\begin{array}{rcl} D^k(M)=Sec (\wedge ^k (T^1M) )&\rightarrow
&A^k(M)\oplus A^{k-1}(M) \\ D&\mapsto & D^0\oplus D^1\cong
D^0+I_M\wedge D^1 ,
\end{array}
\]

where $D^1=i_{\phi _M}D$ and $D^0=D-I_M\wedge D^1$.

As for $A(M)$, we can define on $D(M)=\oplus _{k\in \Z}D^k(M)$ a canonical
Schouten-Jacobi bracket $\lcf \cdot ,\cdot  \rcf ^1 _M:D^k(M)\times
D^r(M)\to D^{k+r-1}(M)$ (see \cite{GM,IM})
\begin{equation}\label{first-order}
\begin{array}{rcl}
\lcf P^0+I_M\wedge P^1,Q^0+I_M\wedge Q^1\rcf _M^1&=&\\[5pt] & &
\kern-82pt \lcf P^0, Q^0\rcf _M + (k-1)P^0\wedge
Q^1+(-1)^k(r-1)P^1\wedge Q^0
\\ &&\kern-82pt +I_M\wedge \Big ( \lcf P^1,Q^0\rcf _M -(-1)^k \lcf P^0,
Q^1\rcf _M +(k-r)P^1\wedge Q^1 \Big ),
\end{array}
\end{equation}
for $P=P^0+I_M\wedge P^1\in D^k(M)$ and $Q=Q^0+I_M\wedge Q^1\in
D^r(M)$. The bracket $\lcf \cdot , \cdot \rcf _M^1$ is the unique
graded bracket characterized by:

{\it i)} It extends the Lie bracket on $D^1(M)$ defined by
(\ref{firstorder});

{\it ii)} $\lcf X\oplus f,g \rcf ^1_M =X(g)+fg$, for $X\oplus f
\in D^1(M)$ and $g\in C^\infty (M)$;

{\it iii)} $\lcf D,E\rcf ^1 _M =-(-1)^{(p-1)(q-1)}\lcf E,D\rcf _M$, for
$D\in A^{p}(M)$, $E\in A^{q}(M)$.

{\it iv)} $\lcf D,E\wedge F\rcf ^1 =\lcf D,E\rcf ^1_M\wedge
F+(-1)^{(p-1)q}E\wedge\lcf D,F\rcf ^1 _M-(i_{\phi _M}D)\wedge E\wedge F$,

for $D\in D^{p}(M)$, $E\in D^{q}(M)$ and $F\in D^{\ast }(M)$;

{\it v)} $(-1)^{(p-1)(r-1)}\lcf\lcf D\kern-0.2pt,\kern-0.2ptE\rcf
^1_M,\kern-0.2pt F\rcf ^1_M+(-1)^{(p-1)(q-1)}\lcf\lcf
E\kern-0.2pt,\kern-0.2pt F \rcf ^1_M,\kern-0.2pt D \rcf
^1_M+(-1)^{(q-1)(r-1)}$ $\lcf \lcf F\kern-0.2pt,\kern-0.2ptD\rcf
^1_M,$ $\kern-0.2pt E\rcf ^1_M=0$, for $D\in D^{p}(M)$, $E\in
D^{q}(M)$ and $F\in D^{r}(M)$.

On the other hand, the space $\Theta^k(M)=Sec(\wedge^k(T^1M)^*)$ of
sections of the vector bundle $\wedge^k(T^1M)^*\to M$ can be identified
with $\Omega^k(M)\oplus \Omega^{k-1}(M)$. Actually, there exists an
isomorphism between $\Theta^k(M)$ and $\Omega^k(M)\oplus \Omega^{k-1}(M)$
given by the formula
\[
\begin{array}{rll}
\Theta^k(M)=Sec(\wedge^k(T^1M)^*)&\to& \Omega^k(M)\oplus
\Omega^{k-1}(M)\\\alpha&\to &\alpha^0\oplus \alpha^1\cong
\alpha^0+\phi_M\wedge \alpha^1
\end{array}
\]
where
\[
\alpha^1=i_{I_M}\alpha,\;\;\;\; \alpha^0=\alpha -\phi_M\wedge
\alpha^1.
\]
In other words,
\[
\begin{array}{lcl}
\alpha(X_1\oplus f_1,\dots ,X_k\oplus f_k)&=&\alpha^0(X_1,\dots
,X_k) + \displaystyle\sum_{i=1}^k(-1)^{i+1}f_i\alpha^1(X_1,\dots,
\hat{X}_i,\dots ,X_k)
\end{array}
\]
for $X_1\oplus f_1,\dots ,X_k\oplus f_k\in Sec (T^1M).$

As for $\Omega(M),$ we can define on $\Theta(M)=\oplus_{k\in
\Z}\Theta^k(M)$ the Jacobi differential $d_M^1:\Theta^k(M)\to
\Theta^{k+1}(M)$ as the map characterized by the following
properties:
\begin{enumerate}
\item $d^1_M$ is a $\R$-linear map.
\item If $f\in C^\infty(M)$ and $j^1f\in Sec((T^1M)^*)$ is the
first jet prolongation of $f$ then $d_M^1f=j^1f$.
\item $d^1_M(\alpha\wedge \beta)=d^1_M\alpha\wedge \beta +
(-1)^p\alpha\wedge d^1_M\beta-\phi_M\wedge \alpha\wedge \beta$,
for $\alpha\in \Theta^p(M)$ and $\beta\in \Theta^q(M)$.
\item $(d^1_M)^2=0,$ that is, $d^1_M$ is a cohomology operator.
\end{enumerate}
Under the isomorphism between $\Theta^k(M)$ and $\Omega^k(M)\oplus
\Omega^{k-1}(M)$ the operator $d_M^1$  is given by
\[d_M^1(\alpha^0,\alpha^1)=(d_M\alpha^0,-d_M\alpha^1+ \alpha^0),
\]
for $(\alpha^0,\alpha^1)\in \Omega^k(M)\oplus \Omega^{k-1}(M)\cong
\Theta^k(M).$

 To finish
with this section, we recall that it is easy to identify $P\in
A^k(M)$ (resp., $D=D^0+I_M\wedge D^1\in D^k(M)$) with a
polyderivation $\{ \cdot,\ldots, \cdot \} _P \colon C^\infty
(M)\times \stackrel{k)}{\ldots}\times C^\infty (M)\to C^\infty
(M)$ (resp., a first-order polydifferential operator $\{
\cdot,\ldots, \cdot \} _D\colon C^\infty (M)\times
\stackrel{k)}{\ldots}\times C^\infty (M)\to C^\infty (M)$) given
by

\begin{equation}\label{bracket1}
\{ f_1,\ldots, f_k \} _P=\langle P, df_1\wedge \ldots \wedge df_k
\rangle
\end{equation}
(resp.,
\begin{equation}\label{bracket2}
\begin{array}{rcl}
\{ f_1,\ldots, f_k \} _D&=&\langle D, j^1f_1\wedge \ldots \wedge
j^1f_k \rangle =\{ f_1, \ldots ,f_k \}_{D^0}
\\ &&+ \displaystyle \sum _{i=1}^k (-1)^{i+1}f_i\,\{
f_1,\ldots ,\widehat{f_i}, \ldots ,f_k \}_{D^1})
\end{array}
\end{equation}
for all $f_1,\ldots ,f_k\in C^\infty (M)$.  Note that $(A(M),\lcf \, ,\,
\rcf _M)$ is naturally embedded into $(D(M),\lcf \, ,\, \rcf ^1_M)$.
Actually, elements of $(A(M),\lcf \, ,\, \rcf _M)$ are just those
$D\in(D(M),\lcf \, ,\, \rcf ^1_M)$ for which $i_{\phi _M}D=0$.

\section{Homogeneous structures}

\setcounter{equation}{0}

\subsection{Homogeneous tensors}

In this Section we will consider a particular class of tensors
related to a distinguished vector field on a manifold.

Let $M$ be a differentiable manifold and let $\Delta$ be a vector field on
$M$. The pair $(M,\Delta )$ will be called a {\em homogeneous structure}.

A function $f\in C^\infty (M)$ is {\em $\Delta$-homogeneous of degree
$n$}, $n\in \R$, if $\Delta (f)=n\, f$. The space of $\Delta$-homogeneous
functions of degree $n$ will be denoted by $S^n_\Delta (M)$. Similarly, a
tensor $T$ is {\em $\Delta$-homogeneous of degree $n$} if $\Li_\Delta T =
nT$. Here $\Li$ denotes the Lie derivative. In particular, $\Delta$ itself
is homogeneous of degree zero. As a result of properties  of the Lie
derivative we get the following properties of the introduced homogeneity
gradation.

    \begin{enumerate}

\item The tensor product $T\otimes S$ of $\Delta$-homogeneous tensors
of degrees $n$ and $m$ respectively, is homogeneous of degree
$n+m$.

\item The contraction of tensors of homogeneity degrees $n$ and $m$ is
homogeneous of degree $n+m$.

\item The exterior derivative preserves the homogeneity degree of forms.
\item The Schouten-Nijenhuis  bracket of multivector fields of homogeneity
degrees $n$ and $m$ is homogeneous of degree $n+m$.
    \end{enumerate}

These properties justify our choice of the homogeneity gradation, which is
compatible with the polynomial gradation introduced in \cite{TU} and
differs by a shift from  homogeneity gradation of contravariant tensors in
some other papers (e.g \cite{Li}).

\begin{example}

{\rm {\it i)} The simplest example of a homogeneous structure is
the pair $(N\times \R,
\partial _s)$, where $\partial _s$ is the canonical vector field
on $\R$. $(N\times \R, \partial _s)$ will be called a {\em free
homogeneous structure}. In this case,
\[
\begin{array}{rcl}
S^n_\Delta (M)\kern-5pt&\kern-2pt=&\kern-8pt \{ f\in C^\infty
(N\times\R) \colon \, f(x,s)\kern-1pt=\kern-1pt \mbox{e}^{ns}
f_N(x), \mbox{ with }f_N\in C^\infty (N), \, \forall\, (x,s)\in
N\times \R\} .
\end{array}
\]

{\it ii)} Let $M=N\times\R$ and $\Delta=s\partial _s$, $s$ being
the usual coordinate on $\R$. In this case

\[
\begin{array}{rcl}
S^n_\Delta (M)\kern-5pt&\kern-2pt=&\kern-8pt \{ f\in C^\infty
(N\times \R) \colon \, f(x,s)\kern-1pt=\kern-1pt s^n f_N(x),
\mbox{ with }f_N\in C^\infty (N), \, \forall\, (x,s)\in N\times
\R\}.
\end{array}
\]

{\it iii)} If $M=\R$ and $\Delta =s^2\partial _s$, then
$S^0_{\Delta}(M)=\R$ and $S^n_\Delta (M)=\{0\}$ for $n\ne 0$
because the differential equation $s^2\frac{\partial f}{\partial
s}=nf$ has no global smooth solutions on $\R$ for $n\ne 0$. }
\end{example}

Using coordinates adapted to the vector field, one can easily
prove the following result.
\begin{proposition}\label{transverso}
Let $(M,\Delta )$ be a homogeneous structure and $N$ be a closed
submanifold in $M$ of codimension 1 such that $\Delta$ is
transversal to $N$. Then, there is a tubular neighborhood $U$ of
$N$ in $M$ and a diffeomorphism of $U$ onto $N\times\R$ which maps
$\Delta _{|U}$ into $\partial _s$
\end{proposition}

Let us introduce a particular class of homogeneous structures
which will be important in the sequel.
\begin{definition}
A  homogeneous structure $(M,\Delta )$ is said to be {\em strict}
if there is an open-dense subset  $O\subset M$ such that for $x\in
O$
\[
T^\ast _xM=\{ df (x)\colon \, f\in S^1_{\Delta }(M)\}.
\]

\end{definition}
\begin{example}\label{e3.4}
{\rm {\it i)} It is almost trivial that free homogeneous
structures are strict homogeneous.

{\it ii)} An example of a strict homogeneous structure with
$\Delta$ vanishing on a submanifold is the following. Let $E\to M$
be a vector bundle (of rank $>0$) over $M$ and let $\Delta =\Delta
_E$ be the {\em Liouville} vector field on $E$. Then, for
$n\in\Z_+$, $S^n_{\Delta }(E)$ consists of smooth functions on $E$
which are homogeneous polynomials of degree $n$ along fibres. In
particular, functions from $S^1_{\Delta}(E)$ are linear on fibres,
hence generate $T^\ast E$ over $E_0$, the bundle $E$ with the
zero-section removed. }
\end{example}

Now, generalizing the situation for tensors, we will consider
first-order polydifferential operators.

For a homogeneous structure $(M,\Delta )$, we say that $D\in
D^k(M)$ is {\em $\Delta$-homogeneous of degree $n$} if $\lcf
\Delta ,D\rcf ^1_M=nD$. For $P\in A^k(M)$ interpreted as an
element of $D(M)$, it is $\Delta$-homogeneous of degree $n$ when
$\lcf \Delta ,P\rcf _M=\Li_\Delta P=n P$, i.e. the introduced
gradation is compatible with the gradation for tensors. It is easy
to see, using (\ref{first-order}), that $P=P^0+I_M\wedge P^1\in
D^k(M)$ is $\Delta$-homogeneous of degree $n$ if and only if
$P^0\in A^k(M)$  and $P^1\in A^{k-1}(M)$ are $\Delta$-homogeneous
of degree $n$. In particular, the identity operator is homogeneous
of degree zero.

Elements of $D^k(M)$  which are $\Delta$-homogeneous of degree
$1-k$ we will call simply {\em $\Delta$-homogeneous}.

\begin{proposition}

Suppose that $D\in D^k(M)$ is $\Delta$-homogeneous of degree $n$
and suppose that $D'\in D^{k'}(M)$ is $\Delta$-homogeneous of
degree $n'$. Then,
\begin{itemize}
\item[{\it i)}] $D\wedge D'$ is $\Delta$-homogeneous of degree
$n+n'$.
\item[{\it ii)}] $\lcf D, D'\rcf ^1_M$ is $\Delta$-homogeneous of degree
$n+n'$.
\end{itemize}
\end{proposition}

{\bf Proof.-} These properties are immediate consequences of
properties of the Schouten-Jacobi bracket $\lcf \cdot , \cdot \rcf
^1_M$ (see Section \ref{graded}) and the fact that $i _{\phi
_M}\Delta =0$. \hfill$\Box$

We can characterize homogeneous operators for strict homogeneous
structures in terms of the corresponding $k$-ary brackets as
follows.

\begin{proposition}\label{p3.6}
Let $(M,\Delta )$ be a strict homogeneous structure. Then,
\begin{itemize}
\item[{\it i)}] $P\in A^k(M)$ is $\Delta$-homogeneous of degree
$n$ if and only if $\{ f_1,\ldots ,f_k \}_P$ is
$\Delta$-homoge\-neous of degree $n+k$, for all $f_1,\ldots ,f_k
\in
S^1_{\Delta } (M)$, where $\{ \cdot ,\ldots ,\cdot \}_P$ is the
bracket defined as in (\ref{bracket1}).
\item[{\it ii)}] $D\in D^k(M)$ is $\Delta$-homogeneous of degree
$n$ if and only if $\{ f_1,\ldots ,f_k \}_D$ is
$\Delta$-homogeneous of degree $n+deg(f_1)+\ldots +deg(f_k),$ for
all $\Delta$-homogeneous functions $f_1,\ldots ,f_k$ of degree $1$
or $0$.
\end{itemize}
\end{proposition}
{\bf Proof.-} {\it i)} follows from the identity

\[
\Delta (\{ f_1,\ldots ,f_k \}_P )= \langle \lcf{\Delta}, P\rcf_M,
df_1\wedge \ldots \wedge df_k\rangle + \langle P, {\cal L}_{\Delta
}(df_1\wedge \ldots \wedge df_k )\rangle ,
\]
for $f_1,\ldots ,f_k\in C^\infty (M)$, where ${\cal L}$ denotes
the usual Lie derivative operator, and the fact that $df_1\wedge
\ldots \wedge df_k$, with $\Delta$-linear functions $f_1,\ldots
,f_k$, generate $\wedge ^kT^\ast M$ over an open-dense subset.

The proof of {\it ii)} is analogous. \hfill$\Box$

Next, we will consider the particular case when $\Delta$ is the
Liouville vector field $\Delta_E$ on a vector bundle $E$. We
recall that, in such a case, $S_{\Delta_E}^1(E)$ is the space of
linear functions on $E$ and $S_{\Delta_E}^0(E)$ is the space of
basic functions on $E$ (see Example \ref{e3.4}).

\begin{corollary}
Let $E\to M$ be a vector bundle over $M$, $\Delta_E$ be the
Liouville  vector field on $E$  and $(E,\Delta_E)$ be the
corresponding strict homogeneous structure. Then:
\begin{enumerate}
\item $P\in A^k(E)$ is $\Delta_E$-homogeneous if and only if $P$
is linear, that is,
\begin{equation}\label{1}
\{f_1,\dots ,f_k\}_P\in S_{\Delta_E}^1(E), \mbox{ for }f_1,\dots
,f_k\in S^1_{\Delta_E}(E).
\end{equation}
\item $D\in D^k(M)$ is $\Delta_E$-homogeneous if and only if
\begin{equation}\label{2}
\begin{array}{ll}\{f_1,\dots ,f_k\}_D\in S_{\Delta_E}^1(E), &\mbox{ for
}f_1,\dots ,f_k\in S^1_{\Delta_E}(E),\\ \{1,f_2,\dots ,f_k\}_D\in
S_{\Delta_E}^0(E),& \mbox{ for }f_2,\dots ,f_k\in
S^1_{\Delta_E}(E).\end{array}
\end{equation}
\end{enumerate}
\end{corollary}
{\bf Proof.-} $(i)$ follows from Proposition \ref{p3.6}.

On the other hand, if $D\in D^k(M)$ is $\Delta_E$-homogeneous
then, using again Proposition \ref{p3.6}, we deduce that (\ref{2})
holds.

Conversely, suppose that (\ref{2}) holds.

Then, if $f_1^0\in S_{\Delta_E}^0(E)$ and $f_1^1,\dots ,f_k^1\in
S_{\Delta_E}^1(E),$ we have that \[
S^1_{\Delta_E}(E)\ni\{f_1^0f_1^1, f_2^1,\dots
,f_k^1\}_D=f_1^0\{f_1^1, f_2^1,\dots ,f_k^1\}_D+f_1^1\{f_1^0,
f_2^1,\dots ,f_k^1\}_D-f_1^0f_1^1\{1, f_2^1,\dots ,f_k^1\}_D.
\]
This implies that
\[
f_1^1\{f_1^0, f_2^1,\dots ,f_k^1\}_D\in S_{\Delta_E}^1(E), \;\;\;
\forall f_1^1\in S_{\Delta_E}^1(E).
\]
Thus,
\begin{equation}\label{3}
\{f_1^0, f_2^1,\dots ,f_k^1\}_D\in S_{\Delta_E}^0(E),\mbox{ for
}f_1^0\in S_{\Delta_E}^0(E)\mbox{ and } f^1_2,\dots ,f_k^1\in
S_{\Delta_E}^1(E).
\end{equation}
Now, we will  see that
\begin{equation}\label{4}
\{1, f_2^0,f_3^1,\dots ,f_k^1\}_D=0,\mbox{ for }f_2^0\in
S_{\Delta_E}^0(E)\mbox{ and } f^1_3,\dots ,f_k^1\in
S_{\Delta_E}^1(E).
\end{equation}
If $f_2^1\in S_{\Delta_E}^1(E)$, we obtain that
\[ S_{\Delta_E}^0(E)\ni \{1, f_2^0f_2^1,f_3^1,\dots
,f_k^1\}_D=f_2^0\{1, f_2^1,f_3^1,\dots ,f_k^1\}_D+f_2^1\{1,
f_2^0,f_3^1, \dots ,f_k^1\}_D.
\]
Therefore, we deduce that
\[
f_2^1\{1, f_2^0,f_3^1,\dots ,f_k^1\}_D\in S_{\Delta_E}^0(E),
\;\;\forall f_2^1\in S^1_{\Delta_E}(E),
\]
and, consequently,
\[
\{1, f_2^0, f_3^1, \dots ,f_k^1\}_D=0.
\]
Next, we will prove that
\begin{equation}\label{5}
\{f_1^0, f_2^0, f_3^1,\dots ,f_k^1\}_D=0,\mbox{ for }f_1^0,
f_2^0\in S_{\Delta_E}^0(E)\mbox{ and } f^1_3,\dots ,f_k^1\in
S_{\Delta_E}^1(E).
\end{equation}
If $f_2^1\in S_{\Delta_E}^1(E)$ then, using (\ref{3}) and
(\ref{4}), we have that
\[ S_{\Delta_E}^0(E)\ni \{f_1^0, f_2^0f_2^1,f_3^1,\dots
,f_k^1\}_D=f_2^0\{f_1^0, f_2^1,f_3^1, \dots
,f_k^1\}_D+f_2^1\{f_1^0, f_2^0,f_3^1,\dots ,f_k^1\}_D.
\]
This implies that
\[
f_2^1\{f_1^0,f_2^0,f_3^1,\dots ,f_k^1\}_D\in S_{\Delta_E}^0(E),
\;\;\forall f_2^1\in S_{\Delta_E}^1(E),
\]
and thus (\ref{5}) holds.

Proceeding as above, we also  may deduce that
\[
\{f_1^0,\dots ,  f_r^0, f_{r+1}^1,\dots ,f_k^1\}_D=0,
\]
for $f_1^0,\dots ,f_r^0\in S_{\Delta_E}^0(E)$ and $f_{r+1}^1,\dots
,f_k^1\in S_{\Delta_E}^1(E),$ with $2\leq r\leq k.$

Therefore, $D$ is $\Delta_E$-homogeneous (see Proposition
\ref{p3.6}). \hfill $\Box$

\begin{remark}
{\rm We remark that Poisson (Jacobi) structures which are
homogeneous with respect to the Liouville vector field of a vector
bundle play an important role in the study of mechanical systems.
Some examples of these structures are the following: the canonical
symplectic structure on the cotangent bundle $T^*M$ of a manifold
$M$, the Lie-Poisson structure on the dual space of a real Lie
algebra of finite dimension, and the canonical contact structure
on the product manifold $T^*M\times \R$ (for more details, see
\cite{IM0}).}
\end{remark}

\subsection{Poisson-Jacobi reductive structures}

\begin{definition}
A {\em Poisson-Jacobi (PJ) reductive structure} is a triple
$(M,N,\Delta )$, where $(M,\Delta )$ is a homogeneous structure
and $N$ is a 1-codimensional closed submanifold of $M$ such that
$\Delta$ is transversal to $N$.
\end{definition}
 From Proposition \ref{transverso}, we deduce the following result.
\begin{proposition}\label{p3.8}
Let $(M,N,\Delta )$ be a PJ reductive structure. Then, there is  a
tubular neighborhood $U$ of $N$ in $M$ such that $(U,N,\Delta
_{|U})$ is diffeomorphically equivalent to the free PJ reductive
structure $(N\times \R,N,\partial _s)$.
\end{proposition}
Now, we pass to the main result of the paper.

Let $(M,N,\Delta)$ be a PJ reductive structure. Let us consider a tubular
neighborhood $U$ of $N$, like in Proposition 3.8. There is the unique
function $\tilde{1}_N\in S^1_{\Delta}(U)$ such that
$(\tilde{1}_N)_{|N}\equiv 1$. Under the diffeomorphism between $U$ and
$N\times \R$, $\tilde{1}_N$ is the positive function on $N\times \R$ \[
(x,s)\in N\times \R\to e^s\in \R.\] Let us denote by ${\cal F}$ the
foliation defined as the level sets of this function and by $A({\cal F})$,
$D({\cal F})$ the spaces of elements of $A(U)$, $D(U)$ which are tangent
to ${\cal F}$. Here we call $P\in A^k(U)$ {\em tangent to} $\cal F$ if
$P_x\in\wedge^kT_x{\cal F}_x$, where ${\cal F}_x$ is the leaf of $\cal F$
containing $x\in U$. Consequently, $P^0+I_U\wedge P^1\in D^k(U)$ is
tangent to $\cal F$ if $P^0\in A^k(U)$ and $P^1\in A^{k-1}(U)$ are tangent
to $\cal F$.

It is obvious that any $P\in A^k(U)$ has a unique decomposition
$P=P^0_{\cal F}+\Delta _{|U}\wedge P^1_{\cal F}$, where $P^0_{\cal F}\in
A^k(\cal F)$ and $P^1_{\cal F}\in A^{k-1}(\cal F)$. We can use this
decomposition to define, for each $P\in A^k(U)$, operators $J(P)\in
D^k(U)$ and $J_N(P)\in D^k(N)$ by the formulae

    \[ J(P) = P^0_{\cal F} + I_U\wedge P^1_{\cal F}  \]
    and
    \[J_N(P) =J(P)_{|N}. \]

\begin{theorem}\label{main}
Let $(M,N,\Delta )$ be a PJ reductive structure and let $U$ be a
tubular neighborhood of $N$ in $M$ as in Proposition \ref{p3.8}.
Then:
\begin{itemize}

\item[{\it i)}] The mapping $J$  defines a
one-to-one correspondence between $\Delta _{|U}$-homogeneous
multivector fields on $U$ and $\Delta _{|U}$-homogeneous
first-order polydifferential operators on $U$ which are tangent to
the foliation ${\cal F}$;

\item[{\it ii)}]  The mapping $J_N$  defines a
one-to-one correspondence between $\Delta _{|U}$-homogeneous
multivector fields on $U$ and first-order polydifferential
operators on $N$.
\end{itemize}

Moreover,
\begin{enumerate}
\item[$(a)$]
$\{ f_1,\ldots ,f_k \}_P =\{ f_1,\ldots ,f_k \}_{J(P)}$ and $(\{
f_1,\ldots ,f_k \}_P)_{|N} =\{ f_1{}_{|N},\ldots ,f_k{}_{|N}
\}_{J_N(P)}$
\item[$(b)$]
$\lcf J(P),J(Q)\rcf ^1_U=J(\lcf P,Q \rcf _U)$ and  $\lcf
J_N(P),J_N(Q)\rcf^1_N=J_N(\lcf P,Q \rcf _U),$
\end{enumerate}
for all $f_1,\ldots ,f_k\in S^1_\Delta (U)$ and $\Delta
_{|U}$-homogeneous tensors $P,Q\in A(U)$.
\end{theorem}

{\bf Proof.-} The tensors $J(P)$ and $J_N(P)$ clearly satisfy
$(a)$.

Note that the foliation ${\cal F}$ is $\Delta$-invariant, since
$\tilde{1}_N$ is $\Delta$-homogeneous. This implies that $\lcf
\Delta ,A({\cal F})\rcf _M \subset A({\cal F})$, so that $\lcf
\Delta_{\mid U} ,P^0_{\cal F}\rcf _U +\Delta _{|U} \wedge \lcf
\Delta_{\mid U} ,P^1_{\cal F}\rcf _U$ is the decomposition of
$\lcf \Delta_{\mid U} ,P\rcf _U$ for each  tensor $P=P^0_{\cal
F}+\Delta _{|U} \wedge P^1_{\cal F}\in A^k(U)$. This means that if
$P$ is $\Delta_{|U}$-homogeneous then $J(P)$ is also
$\Delta_{|U}$-homogeneous. Conversely, for a  pair $P^0\in A^k
({\cal F}),\ P^1\in A^{k-1}({\cal F})$, $\Delta_{|U}$-homogeneous
of degree $1-k$, the operator $P=P^0 + \Delta_{|U}\wedge P^1$ is
$\Delta_{|U}$-homogeneous. Thus, $J$ is bijective.

Now, due to the fact that for homogeneous $P$, $\lcf \Delta_{\mid U}
,P\rcf _U=(1-k)P=\lcf I_U ,P\rcf ^1_U$, we get by direct calculations
using the properties of the Schouten-Jacobi bracket that $(b)$ is
satisfied.

To prove (ii) we notice first that for a $\Delta_{|U}$-homogeneous
$P$, the operator $(\tilde{1}_N)^{k-1} J(P)$ is homogeneous of
degree zero, i.e. it is $\Delta_{|U}$-invariant. It follows that
$(\tilde{1}_N)^{k-1} J(P)$ and $J(P)$ are uniquely determined by
$J_N(P)$. To show that $J_N$ is surjective, let us take
$D_N=P^0_N+I_N\wedge P^1_N\in D^k(N)$. There are unique
$\bar{P}^0\in A^k(U)$, $\bar{P}^1\in A^{k-1}(U)$ which are
$\Delta_{|U}$-invariant and  equal to $P^0_N$ and $P^1_N$,
respectively, when restricted to  $N$. We just use the flow of
$\Delta_{\mid U}$ to extend tensors on $N$ to $\Delta_{\mid
U}$-invariant tensors on $U$. Then
$\tilde{P}^0=(\tilde{1}_N)^{1-k}\bar{P}^0$ and
$\tilde{P}^1=(\tilde{1}_N)^{1-k}\bar{P}^1$ give rise to a
$\Delta_{|U}$-homogeneous tensor
$\tilde{P}=\tilde{P}^0+\Delta_{\mid U} \wedge \tilde{P}^1$, with
$J_N(\tilde{P})=D_N$. \hfill$\Box$

\begin{remark}\label{Poissonization}

{\rm {\it i)} The above result is a generalization of the main
theorem in \cite{DLM} which states that $\Delta$-homogeneous
Poisson tensors on $M$ can be reduced to Jacobi structures on $N$.
Indeed if $\Lambda$ is Poisson, then $\lcf \Lambda ,\Lambda \rcf
_{|U}=0$, so $\lcf J_N(\Lambda ),J_N(\Lambda )\rcf ^1_N=0$ which
exactly means that $J_N(\Lambda)$ is a Jacobi structure on $N$
(see \cite{GM,IM}). Actually, it is a sort of a
super-Poissonization. Indeed, the Nijenhuis-Schouten bracket
$\lcf\cdot,\cdot\rcf_M$ on $M$ is a graded (or super) Poisson
bracket, while the Schouten-Jacobi bracket
$\lcf\cdot,\cdot\rcf^1_M$ on $N$ is a graded (or super) Jacobi
bracket (cf. \cite{GM2}).

{\it ii)} We call this construction a Poisson-Jacobi reduction,
since it is a half way of the Poisson-Poisson reduction in the
case when $\Gamma =i_{\phi _N}J_N(\Lambda )$ is the vector field
on $N$ whose orbits have a manifold structure. Then, the bracket
$\{ \cdot ,\ldots ,\cdot \}_{J_N(\Lambda )}$ restricted to
functions which are constant on orbits of $\Gamma$ gives a Poisson
bracket on $N/\Gamma $. In the case when $M$ is symplectic, the
Poisson structure on $N/\Gamma $ obtained in this way is the
standard symplectic reduction of the Poisson structure associated
with a symplectic form $\Omega$ on $M$ with respect to the
coisotropic submanifold $N$. An explicit example of the above
construction is the following one. Suppose that the manifold $M$
is  $\R^{2n}$, the submanifold $N$ is the unit sphere $S^{2n-1}$
in $\R^{2n}$ and the vector field $\Delta$ on $\R^{2n}$ is
\[\Delta=\frac{1}{2}\sum_{i=1}^n(q^i\partial_{q^i} + p_i\partial_{p_i}),
\] where $(q^i,p_i)_{i=1,\dots ,n}$ are the usual coordinates on $\R^{2n}$.
 It is clear that $\Delta$ is transversal to $N$. Actually, the map
  \[\R^{2n}-\{0\}\rightarrow S^{2n-1}\times \R,\;\;\;\;\;\; x\to
  (\frac{x}{\|x\|}, \ln{\|x\|^2})\] is a diffeomorphism of $\R^{2n}-\{0\}$
onto $S^{2n-1}\times \R=N\times \R$ which maps
$\Delta_{|\R^{2n}-\{0\}}$ into $\partial_s$. Thus, we will take as
a tubular neighborhood of $N=S^{2n-1}$ in $M=\R^{2n}$ the open
subset $U=\R^{2n}-\{0\}$. Now, let $\Lambda$ be the $2$-vector on
$M$ defined by \[ \Lambda=\sum_{i=1}^n(\partial_{q^i}\wedge
\partial_{p_i}).\]
$\Lambda$ is the Poisson structure associated with the canonical
symplectic $2$-form $\omega$ on $M=\R^{2n}$ given by
$\omega=\displaystyle\sum_{i=1}^ndq^i\wedge dp_i.$

A direct computation proves that $\Lambda_{|U}$ is a
$\Delta_{|U}$-homogeneous Poisson structure. Therefore, it induces a
Jacobi structure $J_N(\Lambda_{|U})$ on $N=S^{2n-1}$. Note that
$J_N(\Lambda_{|U})$ is just the Jacobi structure associated with the
canonical contact $1$-form $\eta$ on $S^{2n-1}$ defined by
\[
\eta=\frac{1}{2}j^*(\displaystyle\sum_{i=1}^{n}(q^idp_i-p_idq^i)),
\]
where $j:S^{2n-1}\to \R^{2n}$ is  the canonical inclusion (for the
definition of the Jacobi structure associated with a contact
$1$-form, see, for instance,  \cite{ChLM}). This Poisson-Jacobi
reduction can be associated also with a reduction with respect to
a Hamiltonian action of $S^1$ on $\R^{2n}$. Indeed, consider the
harmonic oscillator Hamiltonian $H:\R^{2n}\to \R$ given by
\[H=\frac{1}{2}\displaystyle\sum_{i=1}^{n}((q^i)^2 + (p_i)^2)\]
and the Hamiltonian vector field ${\cal
H}_H^\Lambda=i_{dH}(\Lambda)$ of $H$ with respect to $\Lambda$,
that is,
\[ {\cal H}_H^{\Lambda}=\sum_{i=1}^n(p_i{\partial_{q^i}}-q^i{\partial_{p_i}}). \]
The orbit of ${\cal H}_H^\Lambda$
passing through $(q^i,p_i)$ is the curve $\alpha_{(q^i,p_i)}:
\R\to \R^{2n}$ on $\R^{2n}$
\[
\begin{array}{lcl}\alpha_{(q^i,p_i)}(t)&=&(q^1\cos{t} + p_1\sin{t},\dots,
q^n\cos{t} + p_n\sin{t},
\\&&p_1\cos{t}-q^1\sin{t},\dots ,  p_n\cos{t}-q^n\sin{t}).
\end{array}
\]
Consequently, $\alpha_{(q^i,p_i)}$ is  periodic with period $2\pi$ which
implies that the flow of ${\cal H}_H^\Lambda$ defines a symplectic action
of $S^1$ on $\R^{2n}$ with the momentum map given by $H$. Moreover, the
restriction $\Gamma$ of ${\cal H}_H^\Lambda$ to $S^{2n-1}$ is tangent to
$S^{2n-1}$ and $\Gamma$ is a regular vector field on $S^{2n-1}$, that is,
the space of orbits of $\Gamma$, $S^{2n-1}/\Gamma$, has a manifold
structure and, thus, $S^{2n-1}/\Gamma\cong S^{2n-1}/S^1$ is a symplectic
manifold. Actually, the reduced symplectic space $S^{2n-1}/S^1$ is the
complex projective space with the standard symplectic structure.

{\it iii)} We call the inverse of the map $P\mapsto J_N(P)=D_N$
the {\em Poissonization} of $D_N\in D^k(N)$. This map is a
homomorphism of the Schouten-Jacobi bracket on $D(N)$ into the
Schouten-Nijenhuis bracket of $\Delta$-homogeneous multivector
fields in a neighborhood of $N$ in $M$. In particular, it maps
Jacobi structures into Poisson structures. For free PJ reductive
structures we get, like in \cite{DLM} for the case $k=2$, that the
Poissonization of $D_N=P^0_N+I_N\wedge P^1_N$ is
$\mbox{e}^{(1-k)s}(P^0_N+\partial _s\wedge P^1_N)$ on $N\times\R$.
 }

\end{remark}

Using Theorem \ref{main} and generalizing Remark
\ref{Poissonization} {\it i)}, we have the following result which
relates homogeneous Nambu-Poisson tensors on $M$ to Nambu-Jacobi
tensors on $N$ (see  \cite{MVV,T} for the definition of a
Nambu-Poisson and a Nambu-Jacobi tensor).
\begin{corollary}

Let $(M,N,\Delta )$ be a PJ reductive structure. For a tubular
neighborhood $U$ of $N$ in $M$ there is a one-to-one
correspondence between $\Delta_{|U}$-homogeneous Nambu-Poisson
tensors on $M$ and Nambu-Jacobi tensors on $N$.
\end{corollary}

{\bf Proof.-} We know that a tensor $P\in A^{k}(M)$ on a manifold
$M$ is Nambu-Poisson if and only if
\begin{equation}\label{Nambu-Poisson}
\lcf\lcf\dots \lcf \lcf P,f_1\rcf _M,f_2\rcf _M,\ldots
,f_{k-1}\rcf _M,P\rcf _M=0 ,
\end{equation}
for $f_1,\ldots f_{k-1}\in C^\infty (M)$ and that $D\in D^{k}(M)$
is a Nambu-Jacobi structure on $M$ if and only if
\begin{equation}\label{Nambu-Jacobi}
\lcf \lcf\dots \lcf \lcf D,f_1\rcf ^1_M,f_2\rcf ^1_M,\ldots
,f_{k-1}\rcf ^1_M ,D\rcf ^1_M=0 ,
\end{equation}
for $f_1,\ldots f_{k-1}\in C^\infty (M)$.

Therefore, our result follows from (\ref{Nambu-Poisson}),
(\ref{Nambu-Jacobi}) and Theorem \ref{main}.\hfill$\Box$

The above result is local. We can get global results in particular
classes. The following one has been proved in \cite{GIMPU} for
bivector fields by a different method.
\begin{theorem}\label{global}
Let $E\to M$ be a vector bundle of rank $n$, $n>$1, and let $A$ be
an affine hyperbundle of $E$, i.e. an affine subbundle of rank
$(n-1)$ and not intersecting the $0$-section of $E$. Then, the
association $P\mapsto J_A(P)$ establishes a one-to-one
correspondence between $\Delta _E$-homogeneous tensors $P\in
A^k(E)$, the vector field $\Delta _E$ being the Liouville vector
field, and those $D_A\in D^k (A)$ which are affine, i.e. such that
$\{ h_1,\dots,h_k\}_{D_A}$ is affine whenever $h_1,\dots,h_k$ are
affine (along fibers) functions on $A$. Moreover, for this
correspondence, \be\label{af1} \lcf J_A(P),J_A(Q)\rcf
^1_A=J_A(\lcf P,Q\rcf _E). \ee
\end{theorem}

{\bf Proof of Theorem \ref{global}.-} The Liouville vector field
$\Delta_E$ is clearly transversal to $A$, so the association
$P\mapsto J_A(P)$ satisfies $$(\{ f_1,\ldots ,f_k\} _P)_{|A}=\{
f_1{}_{|A},\ldots ,f_k{}_{|A}\} _{J_A(P)} $$ and (\ref{af1})
according to Theorem \ref{main}. The affine functions on $A$ are
exactly restrictions of linear functions on $E$ (see the next
Lemma \ref{lema}), so $J_A(P)$ is affine.

Conversely, according to Theorem \ref{main}, there is a neighbourhood U of
$A$ in $E$ on which $\Delta_E$ nowhere vanishes and a $(\Delta_E)_{\mid
U}$-homogeneous $k$-vector field $P_U$ on $U$ such that $D_A=J_A(P_U)$. We
will show that $P_U$ is {\em linear}, i.e. that $\{ (f_1)_{\mid
U},\dots,(f_k)_{\mid U}\}_{P_U}$ is the restriction to $U$ of a linear
function on $E$ for all linear functions $f_1,\dots,f_k$ on $E$. In the
case of a $0$-tensor, i.e. a function $f\in C^\infty(U)$, this means that
$f$ is the restriction to $U$ of a linear function on $E$.

Indeed, since by Theorem \ref{main} $$(\{ (f_1)_{\mid
U},\dots,(f_k)_{\mid U}\}_{P_U})_{|A}=\{ f_1{}_{|A},\ldots
,f_k{}_{|A}\} _{D_A}, $$ the function $\{ (f_1)_{\mid
U},\dots,(f_k)_{\mid U}\}_{P_U}$ is $\Delta_E$-homogeneous on $U$
and its restriction to $A$ is affine, thus it is the restriction
to $U$ of a linear function. Note that every affine function on
$A$ has a unique extension to a linear function on the whole $E$
(see the next Lemma \ref{lema}). Moreover, two
$\Delta_E$-homogeneous functions $f$ and $g$ on $U$ which coincide
on $A$ must coincide on the $\Delta_E$ orbits of points from $A$
and , since $A$ is an affine hyperbundle of $E$ not intersecting
the $0$-section of $E$, we deduce that $f=g$ on $U$.

What remains to be proven is that $P_U$ has a unique extension to a
$\Delta_E$-homogeneous tensor on $E$ that follows from the next Lemma
\ref{extension}.\hfill$\Box$

\begin{lemma}\label{lema}
Let $E$ be a real vector bundle over $M$ and $A$ be an affine
hyperbundle  of $E$ not intersecting the $0$-section $0:M\to E$ of
$E$. Suppose that $A^+$  is the real vector bundle over $M$ whose
fiber at the point $x\in M$ is the real vector space $
A_x^+=Aff(A_x,\R),$ that is, $A_x^+$ is the space of real affine
functions  on $A_x$. Then, the map $R_A:E^*\to A^+$ defined by
$R_A(\alpha_x)=(\alpha_x)_{|A_x},$ for $\alpha_x\in E^*_x$ is an
isomorphism of vector bundles.
\end{lemma}

{\bf Proof.-} Let $x$ be a point of $M$ and $\alpha_x\in E_x^*.$
Then, it is easy to prove that $R_A(\alpha_x)\in A_x^+$ and that
the map $(R_A)_{|E_x^*}:E_x^{*}\to A_x^+$ is linear. Moreover, if
$R_A(\alpha_x)=0$, we have that $(\alpha_x)_{|A_x}=0$ and,  using
that $0(x)\notin A_x$, we conclude that $\alpha_x=0$. Thus,
$(R_A)_{|E_x^*}$ is injective and, since $\dim E_x^*=\dim
A_x^+=n$, we conclude that $(R_A)_{|E_x^*}:E_x^*\to A_x^+$ is a
linear isomorphism. This proves the result. \hfill$\Box$

\begin{lemma}\label{extension}
Let $\tau:E\to M$ be a vector bundle of rank $n$, $n>$1, $A$ be an
affine hyperbundle of $E$ not intersecting the $0$-section of $E$
and $U$  be a neighborhood of $A$ in  $E$. If  $P$ is a
linear-homogeneous $k$-contravariant tensor field on $U$ then $P$
has a unique extension to a $\Delta _E$-homogeneous (linear)
$k$-contravariant tensor field $\tilde{P}$ on $E$.
\end{lemma}

{\bf Proof.-} The statement is local in $M$, so let us choose
local coordinates $x=(x^a)$ in $V\subset M$ and the adapted linear
coordinates $(x^a,\xi_i)$ on $E_{\mid V}$, associated with a
choice of a basis of local sections of $E_{\mid V}$. In these
coordinates, the tensor $P$ can be written in the form
\bea\label{prol}
&P=\sum_{i_1,\dots,i_k}f^k_{\xi_{i_1},\dots,\xi_{i_k}}(x,\xi)\pa_{\xi_{i_1}}
\ot\cdots\ot\pa_{\xi_{i_k}}+\\
&+\sum_{i_1,\dots,i_{k-1},a}f^{k-1}_{\xi_{i_1},\dots,\xi_{i_{k-1}},x^a}(x,\xi)
\pa_{\xi_{i_1}}\ot\cdots\ot\pa_{\xi_{i_{k-1}}}\ot\pa_{x^a}+\nn\\
&+\sum_{i_1,\dots,i_{k-1},a}f^{k-1}_{\xi_{i_1},\dots,\xi_{i_{k-2}},x^a,\xi_{i_{k-1}}}(x,\xi)
\pa_{\xi_{i_1}}\ot\cdots\ot\pa_{\xi_{k-2}}\ot\pa_{x^a}\ot\pa_{\xi_{i_{k-1}}}+
\cdots +\nn\\
&+\sum_{a_1,\dots,a_k}f^0_{x^{a_1},\dots,x^{a_k}}(x,\xi)\pa_{x^{a_1}}\ot\cdots
\ot\pa_{x^{a_k}}.\nn \eea By linearity of the tensor $P$,
$\{\xi_{i_1},\dots,\xi_{i_k}\}_P=f^k_{\xi_{i_1},\dots,\xi_{i_k}}(x,\xi)$
is linear in $\xi$, so it can be extended uniquely to a linear
function on the whole $E_{\mid V}$. Similarly, proceeding by
induction with respect to $m$ one can show that the linearity of
$$\{\xi_{i_1},\dots,x^{a_1}\cdot\xi_{j_1},\dots,x^{a_m}\cdot\xi_{j_m},
\dots,\xi_{i_{k-m}}\}_P $$ implies that
\be\label{lem}f^{k-m}_{\xi_{i_1},\dots,x^{a_1},\dots,x^{a_m},
\dots,\xi_{i_{k-m}}}(x,\xi)\cdot\xi_{j_1}\cdots\xi_{j_m} \ee is
linear for all $j_1,\dots,j_m$. Once we know that (\ref{lem}) are
linear, it is easy to see that

\be\label{lem1}f^{k-1}_{\xi_{i_1},\dots,x^{a_1},
\dots,\xi_{i_{k-1}}}(x,\xi) \ee is constant on fibers, so it
extends uniquely to a function which is constant on the fibers of
$E_{|V}$. On the other hand, since $n>1$ and $U$ is a neighborhood
of $A$ in $E$, there exist $i_1,\dots ,i_{n-1}\in \{1,\dots ,n\}$
such that
\[
U\cap \{\xi_{i_k}=0\}\not=\emptyset, \;\;\; \forall k\in \{1,\dots
,n-1\}.
\]
Using this fact and the linearity of (\ref{lem}), we deduce that
\[f^{k-m}_{\xi_{i_1},\dots,x^{a_1},\dots,x^{a_m},
\dots,\xi_{i_{k-m}}}(x,\xi)=0,\;\;\;\; \mbox{for $m>1$}.
\]
Note that if $rank{(E)}=1$, we have that $\xi_{i_l}=\xi$ and there
 is another possibility, namely
$$f^{k-m}_{\xi_{i_1},\dots,x^{a_1},\dots,x^{a_m},
\dots,\xi_{i_{k-m}}}(x,\xi)=g(x)\xi^{1-m}, $$ which clearly does
not prolong onto $E_{\mid V}$ analytically along fibers. Now we
define the prolongation $\tilde{P}_V$ of $P$ to $E_{\mid V}$ by
the formula (\ref{prol}) but with the prolonged coefficients. It
is obvious that this constructed prolongation $\tilde{P}_V$ of $P$
to $E_{\mid V}$ is homogeneous. By uniqueness of this homogeneous
prolongation on every $E_{\mid V}$ for $V$ running through an open
covering of $M$, we get a unique homogeneous prolongation to the
whole $E$. \hfill$\Box$

\noindent
\begin{remark} {\rm The linearity cannot be replaced by
$\Delta_E$-homogeneity in the above lemma. The simplest
counterexample is just the function $f(x)=\vert x\vert$ which is
$x\partial_x$-homogeneous on $U=\R\setminus \{ 0\}$ but it is not
linear on $U$.}\end{remark}

Finally, we will prove  a dual version of Theorem \ref{main}.

Let $(M,N,\Delta)$ be a PJ reductive structure and let $U$ be a tubular
neighborhood of $N$ in $M$ as in Proposition \ref{p3.8}. The space of
sections of the vector bundle $\wedge^k(T^1U)^*\to U$ (respectively,
$\wedge^k(T^1N)^*\to N)$ is $\Omega^k(U)\oplus \Omega^{k-1}(U)$
(respectively, $\Omega^k(N)\oplus \Omega^{k-1}(N)$) and it is  obvious
that any $\alpha\in \Omega^k(U)$ has a unique decomposition
\be\label{*}
\alpha=\tilde{1}_N(\alpha^0 + d(\ln{\tilde{1}_N})\wedge \alpha^1),
\ee
where $(\alpha^0,\alpha^1)\in \Omega^k(U)\oplus \Omega^{k-1}(U)$ and
\[
i_{\Delta_{|U}}\alpha^0=0,\;\;\; i_{\Delta_{|U}}\alpha^1=0.
\]
Indeed, since $i_{\Delta_{|U}}d(\ln{\tilde{1}_N})=1$, the form $\za^1$ is
defined by $\za^1=(\wt{1}_N)^{-1}i_{\Delta_{|U}}\alpha$ and
$\za^0=(\wt{1}_N)^{-1}\za-d(\ln{\tilde{1}_N})\wedge \alpha^1$. We can use
this decomposition to define, for each $\alpha\in \Omega^k(U)$, a section
$\Psi(\alpha)$ of the vector bundle $\wedge^k(T^1U)^*\to U$ by the formula
\[
\Psi(\alpha)=(\alpha^0,\alpha^1).
\]
On the other hand, a section $(\alpha^0,\alpha^1)\in \Omega^k(U)\oplus
\Omega^{k-1}(U)$ is said to be $\Delta_{|U}$-basic if $\alpha^0$ and
$\alpha^1$ are basic forms with respect to $\Delta_{|U}$, that is,
\[
i_{\Delta_{|U}}\alpha^0=0,\;\;\; i_{\Delta_{|U}}\alpha^1=0,\;\;\; {\cal L}
_{\Delta_{|U}}\alpha^0=0,\;\;\; {\cal L}_{\Delta_{|U}}\alpha^1=0.
\]
In addition, we will denote by $j:N\to U$ the canonical inclusion
and by $\Psi_N:\Omega^k(U)\to \Omega^k(N)\oplus \Omega^{k-1}(N)$
the map defined by
\[
\Psi_N(\alpha)=(\alpha^0_N,\alpha^1_N),\;\;\mbox{ for }\alpha\in
\Omega^k(U),
\]
where $\alpha^0_N=j^*(\alpha)$, $\alpha^1_N=j^*(i_{\Delta_{|U}}\alpha)$.
On the other hand, from (\ref{*}), it follows that
\begin{equation}\label{+}
j^*\alpha=j^*\alpha^0,\;\;\;\;j^*(i_{\Delta_{|U}}\alpha)=j^*\alpha^1,
\end{equation}
(note that $j^*(\tilde{1}_N)$ is the constant function $1$ on $N$), so
$\alpha^0_N=j^*(\alpha^0)$ and $\alpha^1_N=j^*(\alpha^1)$.

\begin{theorem}
Let $(M,N,\Delta)$ be a PJ reductive structure and let $U$ be a tubular
neighborhood of $N$ in $M$ as in Proposition \ref{p3.8}. Then:
\begin{enumerate}
\item
The map $\Psi$ defines a one-to-one correspondence between the space of
$k$-forms on $U$ which are $\Delta_U$-homogeneous of degree $1$ and the
space of sections of the vector bundle $\wedge^k(T^1U)^*\to U$ which are
$\Delta_{|U}$-basic.

\item
The map $\Psi_N$ defines a one-to-one correspondence between the
space of k-forms on $U$ which are $\Delta_{|U}$-homogeneous of
degree $1$ and the space of sections of the vector bundle
$\wedge^k(T^1N)^*\to N,$ that is, $\Omega^k(N)\oplus
\Omega^{k-1}(N).$

Moreover, if $\alpha\in \Omega^k(U)$ is $\Delta_{|U}$-homogeneous of
degree $1$ then
\[
\Psi(d_U\alpha)=d^1_U(\Psi\alpha),\;\;\;
\Psi_N(d_U\alpha)=d_N^1(\Psi_N\alpha),
\]
where $d_U$ is the usual exterior differential on $U$ and $d_U^1$
(respectively, $d_N^1$) is the Jacobi differential on $U$ (respectively,
$N$). \end{enumerate}
\end{theorem}
{\bf Proof.-} Let $\alpha$ be a $k$-form on $U$,
\begin{equation}
\alpha=\tilde{1}_N(\alpha^0 + d(\ln{\tilde{1}_N})\wedge \alpha^1),
\end{equation}
with $(\alpha^0,\alpha^1)\in \Omega^k(U)\oplus \Omega^{k-1}(U)$
satisfying $i_{\Delta_{|U}}\alpha^0=0$ and
$i_{\Delta_{|U}}\alpha^1=0.$ Then
\[
{\cal L}_{\Delta_{|U}}\alpha=\alpha + \tilde{1}_N({\cal
L}_{\Delta_{|U}}\alpha^0 + d(\ln{ \tilde{1}_N})\wedge {\cal
L}_{\Delta_{|U}}\alpha^1).
\]
Thus, since $i_{\Delta_{|U}}({\cal L}_{\Delta_{|U}}\alpha^0)=0$ and
$i_{\Delta_{|U}}({\cal L}_{\Delta_{|U}}\alpha^1)=0$, we conclude that
$\alpha$ is $\Delta_{|U}$-homogeneous of degree $1$ if and only if
$\alpha^0$ and $\alpha^1$ are $\Delta_{|U}$-basic. This proves $(i).$

Since
$$
j^*\alpha=j^*\alpha^0,\;\;\;\;j^*(i_{\Delta_{|U}}\alpha)=j^*\alpha^1,
$$
using $(i)$ and the fact that the map $j^*:\Omega^r(U)\to \Omega^r(N)$
defines a one-to-one correspondence between the space of
$\Delta_{|U}$-basic $r$-forms on $U$ and $\Omega^r(N)$, we deduce $(ii)$.

Finally, if $\alpha\in \Omega^k(U)$ is $\Delta_{|U}$-homogeneous of degree
$1$ then, from (\ref{*}), we obtain that
\[
d_U\alpha=\tilde{1}_N(d_U\alpha^0+d_U(\ln{\tilde{1}_N})\wedge
(\alpha^0-d_U\alpha^1))
\]
and, since
\[
i_{\Delta_{|U}}(d_U\alpha^0)={\cal L}_{\Delta_{|U}}\alpha^0=0,\;\;\;\;\;
i_{\Delta_{|U}}(\alpha^0-d_U\alpha^1)=-{\cal L}_{\Delta_{|U}}\alpha^1=0,
\]
we conclude that (see (\ref{+}))
\[
\Psi(d_U\alpha)=(d_U\alpha^0,\alpha^0-d_U\alpha^1)=d_U^1(\Psi\alpha),
\]
\[
\Psi_N(d_U\alpha)=(d_N(j^*(\alpha^0)),
j^*(\alpha^0)-d_N(j^*(\alpha^1)))=d_N^1(\Psi_N\alpha).
\]
\hfill$\Box$

Using Theorem 3.18, one may recover the following well-known result (see,
for instance, \cite[Proposition 3.58]{MS}).
\begin{corollary} If $\zw$ is a $\Delta_{|U}$-homogeneous of degree 1
symplectic form on $U$, then $\eta=\zw^1_N$ is a contact form on
$N$. The Jacobi structure associated with $\eta$ is $J_N(\zL)$,
where $\zL$ is the $\Delta_{|U}$-homogeneous Poisson tensor
associated with $\zw$.
\end{corollary}
{\bf Proof.-} Since, according to Theorem 3.18,
$$0=\Psi_N(d\zw)=d_N^1(\Psi_N\zw)=(d\zw^0_N,\zw^0_N-d\zw^1_N),$$
we have \be\label{last} d\eta=d\zw^1_N=\zw^0_N=j^*\zw .\ee If the
dimension of $N$ is $2k+1$ then (\ref{last}) implies
$$(d\eta)^{2k}\we\eta =j^*(\zw^{2k}\we i_{\zD_{\mid U}}\zw)=
\frac{1}{k+1}j^*(i_{\zD_{\mid U}}\zw^{2(k+1)}).$$ But
$\zw^{2(k+1)} \neq 0$ on $U$ (the form $\zw$ is symplectic) and
$\zD$ is transversal to $N$, so $j^*(i_{\zD_{\mid U}}\zw^{2(k+1)})
\neq 0$, thus $(d\eta)^{2k}\we\eta \neq 0$ on $N$ and, therefore,
$\eta$ is a contact $1$-form on $N$. The contact form $\eta$
induces an isomorphism of vector bundles $\flat_\eta:TN\ra T^*N$
which on sections takes the form \be\label{last1}
\flat_\eta(X)=\la\eta,X\ran\eta-i_Xd\eta. \ee The Jacobi bracket
$\{ f,g\}_\eta$ induced by $\eta$ is given by $\{ f,g\}_\eta={\cal
H}^\eta_f(g)-g\zG(f)$, where ${\cal H}^\eta_f$ is the `Hamiltonian
vector field' of $f\in C^\infty(N)$ defined by
$$\flat_\eta({\cal H}^\eta_f)=(df-\zG(f)\eta)+f\eta$$
and $\zG$ is the Reeb vector field of $\eta$ determined by
$\flat_\eta(\zG)=\eta$, i.e. $i_\zG d\eta=0$ and
$\la\eta,\zG\ran=1$. Let $\{\cdot,\cdot\}_\zw$ be the Poisson
bracket induced by the symplectic form $\zw$. Due to Theorem 3.11,
it remains to prove that $\{
f,g\}_\eta=(\{\wt{f},\wt{g}\}_\zw)_{\mid N}$, where $\wt{f}$
denotes the unique extension of $f\in C^\infty(N)$ to a $\zD_{\mid
U}$-homogeneous function on $U$. Denote by ${\cal
H}_{\tilde{f}}^\omega$ the Hamiltonian vector field of $\wt{f}$
with respect to $\zw$, i.e. $-i_{{\cal H}^\omega_{\tilde
f}}\zw=d\wt{f}$. It is easy to see that the Reeb vector field of
$\eta$ is $\wt{\zG}_{\mid N}$, $\wt{\zG}={\cal
H}_{\tilde{1}_N}^\omega$, and that ${\cal H}_f^\eta=({\cal
H}_{\tilde f}^\omega+\wt{\zG}(\wt{f})\zD_{\mid U})_{\mid N}$, i.e.
${\cal H}_f^\eta$ is the projection of ${\cal H}^\zw_{\tilde{f}}$
along $\zD$ onto $N$. We have
$$\{ f,g\}_\eta={\cal H}_f^\eta(g)-g\zG(f)=(({\cal H}^\zw_{\tilde{f}}+\wt{\zG}
(\wt{f})\zD_{\mid U})(\wt{g}))_{\mid N}-g\zG(f).
$$
Since $\zD_{\mid U}(\wt{g})=\wt{g}$, it follows that
$$\{f,g\}_\eta= ({\cal H}_{\tilde f}^\zw(\wt{g}))_{\mid N}=
(\{ \wt{f},\wt{g}\}_\zw)_{\mid N}.
$$
\hfill$\Box$
\begin{remark}{\rm If $M=\R^{2n}$, $\Delta$ is the vector field on
$M$ defined by $\Delta=\frac{1}{2}\sum_{i=1}^n (q^i\partial_{q^i}
+ p_i\partial_{p_i}),$ $U$ is the open subset of $M$ given by
$U=\R^{2n}-\{0\},$ $\omega=\sum_{i=1}^n(dq^i\wedge dp_i)$ is the
canonical $\Delta_{|U}$-homogeneous symplectic $2$-form on $U$ and
$N$ is the unit sphere $S^{2n-1}$ in $\R^{2n}$ then $\eta$ is the
canonical contact $1$-form on $S^{2n-1}$ (see Remark
\ref{Poissonization}, $ii)$). }
\end{remark}
\begin{remark}
{\rm A Poisson structure is a particular Lie algebra structure. A useful
generalization of the latter in the graded case is a {\it (strongly)
homotopy Lie algebra} (sh Lie algebra, $L_\infty$-algebra) which appeared
in the works of J.~Stasheff and his collaborators \cite{LM, St0}. Very
close algebraic structures arose in physics  as {\it string products} of
B.~Zwiebach \cite{Zw}. An algebraic background of a homotopy Lie algebra
on a graded vector space $V$ is a graded Lie algebra structure on the
graded space $L(V)=\bigoplus_{n\ge 0}L^n(V)$  of (skew-symmetric)
multilinear maps from $V$ into $V$. The corresponding graded Lie bracket
on $L(V)$ is actually a graded variant of the Nijenhuis-Richardson bracket
$\lcf\cdot,\cdot\rcf^{NR}$ and the homotopy Lie algebra on $V$ is a formal
series $B=\sum_{n\ge 0}B_nh^n$, $B_n\in L^n(V)$ with coefficients which
satisfies the `Master Equation' $\lcf B,B\rcf^{NR}=0$. One requires
additionally that the degree of $B_n$ is $n-2$. Of course, when $B$
reduces to $B_2$, i.e. $B_n=0$ for $n\ne 2$, we deal with a standard
graded Lie bracket on $V$ induced by $B_2:V\ti V\ra V$ of degree 0. When
also $B_1$ is non-trivial, the Jacobi identity for $B_2$ is satisfied only
`up to homotopy'. One can consider this general scheme skipping the
assumption on the degree and one can work with any subalgebra of $L(V)$,
also for non-graded $V$: we just consider the series $B$ with coefficients
in the Lie subalgebra of $L(V)$ and satisfying the Master Equation. Of
course, this general scheme has nothing to do with `homotopy' in general,
when no grading on $V$ or not proper degree of $B_n$ is assumed.

In our case of the Schouten-Nijenhuis and Schouten-Jacobi
brackets, one can consider their homotopy generalizations which
respect the homogeneity, like these brackets do, and obtain the
corresponding Poisson-Jacobi reduction on the level of homotopy
algebras, but the detailed discussion of these problems exceeds
limits of this note and we postpone it to a separate paper.

What we can have for free is the above scheme for the non-graded
case of $V=C^\infty(M)$. The spaces $A^k(M)$ and $D^k(M)$ can be
interpreted as subspaces of $L^n(V)$ and the brackets
$\lcf\cdot,\cdot\rcf_M$ and $\lcf\cdot,\cdot\rcf_M^1$ are
restrictions of $\lcf\cdot,\cdot\rcf^{NR}$ to $A^k(M)$ and
$D^k(M)$, respectively. A {\it formal Poisson} structure on $M$ is
a formal series $B=\sum_{n\ge 0}B_nh^n$, $B_n\in A^n(M)$ such that
$\lcf B,B\rcf_M=0$, where we use the obvious extension of the
Schouten-Nijenhuis bracket to formal series of multivector fields:
$\lcf B,B\rcf_M=\sum_{i,j}\lcf B_i,B_j\rcf_Mh^{i+j-1}$. By
properties of the Schouten-Nijenhuis bracket, only the even part
of $B$ is relevant. If $B_2$ is the only non-trivial part of $B$,
we recognize a standard Poisson structure. If this is the case of
$B_{2k}$, we recognize a {\it generalized Poisson structure} in
the sense of Azc\'arraga, Perelomov, and P\'erez Bueno
\cite{APP1,APP2} (see also \cite{AIP}). Now, according to Theorem
\ref{main}, if $B$ is $\zD$-homogeneous, we can reduce $B$ to a
{\it formal Jacobi} structure on the submanifold $N$ by
$J_N(B)=\sum_{i\ge 0}J_N(B_N)$, since
$$\lcf J_N(B),J_N(B)\rcf_M^1=J_N(\lcf B,B\rcf_M)=0.$$
In particular, this reduces generalized Poisson structures on $M$
to {\it generalized Jacobi structures} on $N$, defined in obvious
way (see \cite{P}). Note also that the corresponding operators
$\pa_B=ad_B$ and $\pa_{J_N(B)}=ad_{J_N(B)}$ act as `homotopy
differentials' in the graded Lie algebras $A^k(M)[[h]]$ and
$D^k(M)[[h]]$, i.e. $\pa_B^2=0$ and $\pa_{J_N(B)}^2=0$,
generalizing the standard Poisson and Jacobi cohomology. }
\end{remark}

\end{document}